\documentclass[11pt]{amsart}

\usepackage{amscd}
\usepackage{amsmath, amssymb}
\usepackage{amsfonts}

\newcommand{\de}{\partial} 
\newcommand{\db}{\overline{\partial}}

\newcommand{\dbar}{\overline{\partial}}
\newcommand{\ddbar}{\sqrt{-1} \partial \overline{\partial}}

\newcommand{\mn}{\sqrt{-1}}

\renewcommand{\leq}{\leqslant}
\renewcommand{\geq}{\geqslant}

\begin{document}
\newcounter{remark}
\newcounter{theor}
\setcounter{remark}{0}
\setcounter{theor}{1}
\newtheorem{claim}{Claim}
\newtheorem{theorem}{Theorem}[section]
\newtheorem{proposition}{Proposition}[section]
\newtheorem{question}{question}[section]
\newtheorem{lemma}{Lemma}[section]
\newtheorem{defn}{Definition}[theor]

\newtheorem{corollary}{Corollary}[section]

\newenvironment{remark}[1][Remark]{\addtocounter{remark}{1} \begin{trivlist}
\item[\hskip
\labelsep {\bfseries #1  \thesection.\theremark}]}{\end{trivlist}}

\title[K-energy on small deformations]{The K-energy on small deformations of constant scalar curvature K\"ahler manifolds}
\author{Valentino Tosatti}
 \address{Department of Mathematics \\ Columbia University \\ New York, NY 10027}
\dedicatory{Dedicated to Professor S.-T. Yau on the occasion of his 60th birthday.}

  \email{tosatti@math.columbia.edu}
\begin{abstract} 
We give a simplified proof of a recent result of X.X. Chen, which together with work of G. Sz\'ekelyhidi implies that on a sufficiently small deformation of a polarized constant scalar curvature K\"ahler manifold the K-energy has a lower bound.
\end{abstract}
\maketitle

\section{Introduction}
The study of canonical metrics in K\"ahler geometry was initiated by Yau \cite{Yau}, and has developed into a very large and active field, see Phong and Sturm \cite{PS} for a survey. A fundamental result in this area says that if a compact K\"ahler manifold $(X,\omega)$ admits a constant scalar curvature K\"ahler (cscK) metric cohomologous to $\omega$, then the Mabuchi K-energy $K_\omega(\varphi)$ of any K\"ahler potential $\varphi$ for $\omega$ is bounded below uniformly. In fact, if the cscK metric is $\omega+\ddbar\varphi_{\mathrm{cscK}}$ then one has
$K_\omega(\varphi)\geq K_\omega(\varphi_{\mathrm{cscK}})$, with equality holding if and only if the metric $\omega+\ddbar\varphi$ is also cscK.

This was first proved for K\"ahler-Einstein metrics by Bando and Mabuchi \cite{BM}. In the groundbreaking papers \cite{D1, D2} Donaldson extended this result to cscK metrics on polarized manifolds (i.e. $[\omega]=c_1(L)$ for some holomorphic line bundle $L$) without nonzero holomorphic vector fields. This was further extended by Chen and Tian \cite{CT2} to cscK metrics on any compact K\"ahler manifold. Later a simpler proof of the lower boundedness of the K-energy for polarized cscK manifolds was provided by Chen and Sun \cite{cs}, and more recently an even shorter proof was found by Li \cite{Chi}.

The main theorem that we want to prove is the following.
\begin{theorem} Let  $(X,J',L')$  be a polarized complex manifold that admits a constant scalar curvature K\"ahler metric in the class $c_1(L')$. If $(X,J,L)$ is a sufficiently small polarized deformation of  $(X,J',L')$, then the K-energy in the class $c_1(L)$ on $(X,J)$ is bounded below.
\end{theorem}

In this situation it follows from Sz\'ekelyhidi's deformation result \cite{Sz1} 
that there exists a smooth test configuration (see section \ref{set} for definitions) with generic fiber $(X,J,L)$ and central fiber with admits a cscK metric in $c_1(L)$ (for a proof, see Proposition 6 in \cite{Sz2}). Then Theorem 1.1 follows immediately from the following:

\begin{theorem}[X.X. Chen \cite{chen}]\label{main}
If $\mathfrak{L}\to \mathfrak{X}\to \mathbf{C}$ is a smooth test configuration with central fiber that admits a cscK metric in $c_1(\mathfrak{L})$, then on the generic fiber the K-energy in the class $c_1(\mathfrak{L})$ is bounded below.
\end{theorem}

In fact, we can compute the infimum of the K-energy in $c_1(L)$ by looking at certain smooth paths of K\"ahler potentials $\varphi_t$ that converge modulo diffeomorphisms to a cscK metric on the central fiber. The infimum of the K-energy on the generic fiber is then equal to the limit of the K-energy along the path $\varphi_t$ when $t$ goes to infinity.

Notice that if the central fiber is not biholomorphic to the generic fiber, then the generic fiber is not K-stable \cite{D3} (since the central fiber has vanishing Futaki invariant \cite{Fu}), and hence it does not admit cscK metrics (at least if it does not have nonzero holomorphic vector fields \cite{D4, stoppa}).

There is an explicit example of such a test configuration, where the central fiber is the Mukai-Umemura threefold \cite{D}, and the generic fiber is Tian's unstable deformation of it \cite{Ti, D5}. This was the first example of a K\"ahler manifold that has a K\"ahler class with K-energy bounded below but without cscK metrics.

The proof of Theorem \ref{main} follows closely the arguments of Theorem 1.7 in \cite{chen}, except that we avoid using a result of Arezzo and Tian \cite{AT} and therefore we do not have to explicitly use weak geodesics in the space of K\"ahler potentials (these are still needed to prove \eqref{chene} below). The key to this simplification is \eqref{important}, which holds for all the paths we consider and not just for geodesics. 

This note is organized as follows: in section \ref{set} we will provide the basic setup concerning smooth test configurations, and in section \ref{prf} we will prove Theorem \ref{main}.\\

\noindent{\bf Acknowledgments. }We are grateful to D.H. Phong, J. Song, J. Sturm, G. Sz\'ekelyhidi and B.Weinkove for very useful discussions and for encouragement, and to S.-T. Yau for his advice and support. This work was partially supported by National Science Foundation grant DMS-1005457.

\section{Setup}\label{set}
A smooth test configuration, as defined by Donaldson in \cite{D3}, is the following data. \begin{itemize}
\item A holomorphic proper submersion
$\pi:\mathfrak{X}\to \mathbf{C}$ with $\mathfrak{X}$ a complex manifold of dimension $n+1$, with a line bundle $\mathfrak{L}\to\mathfrak{X}$ ample on all fibers of $\pi$.
\item An embedding $\mathfrak{X}\subset \mathbf{P}^N\times \mathbf{C}$ with $\mathfrak{L}$ equal to the pullback of the hyperplane bundle and so that this embedding composed with the second projection equals $\pi$.
\item An action of $\mathbf{C}^*$ on $\mathbf{P}^N$ via a $1$-parameter subgroup
$\rho:\mathbf{C}^*\to GL(N+1,\mathbf{C})$, which is extended to an action on $\mathbf{P}^N\times \mathbf{C}$ by acting on the second factor using the product in $\mathbf{C}$, and so that $\mathfrak{X}$ is $\mathbf{C}^*$-invariant and all the maps are $\mathbf{C}^*$-equivariant.
\end{itemize}
In this case, the complex $n$-manifolds $X_\lambda=\pi^{-1}(\lambda)$ with polarization $L_\lambda=\mathfrak{L}|_{X_\lambda}$ are all biholomorphic to a fixed polarized manifold $(X,J,L)$ when $\lambda\neq 0$, while the central fiber $(X_0,J_0,L_0)$ is a complex manifold diffeomorphic to $X$ but usually not biholomorphic to it. Since $c_1(L_0)=c_1(L_\lambda)$ we will identify all these bundles (as complex line bundles), and just call them $L$.

Since we assume that $\pi$ is a submersion, by Ehresmann's theorem we conclude that the family $\mathfrak{X}$ is differentiably trivial, so that there is a diffeomorphism 
$$F:X\times\mathbf{C}\to\mathfrak{X},$$
such that $\pi(F(z,\lambda))=\lambda$. We can think of $F$ as a family of maps
$$F_\lambda:X\to \mathfrak{X},\ \lambda\in\mathbf{C},$$
which are diffeomorphisms with the image $X_\lambda$.
Moreover, from the construction of $F$ in Ehresmann's theorem, we can assume that $F$ extends a given diffeomorphism of the central fiber, and so we may assume that $F_0$ is a biholomorphism between $X$ with the complex structure $J_0$ and its image inside $\mathfrak{X}$.

We can also define a different trivialization of $\mathfrak{X}$ over $\mathbf{C}^*$ using the $1$-parameter subgroup $\rho$. This acts on $\mathfrak{X}\subset \mathbf{P}^N\times \mathbf{C}$ as
$$\rho(\lambda)\cdot(z,\lambda')=(\rho_\lambda(z),\lambda\lambda'),\ \lambda\in\mathbf{C}^*.$$
In particular, if $z$ is in $X_1$ (the fiber over $1$) then $\rho_\lambda(z)$ is in $X_\lambda$, and $\rho_\lambda$ gives a biholomorphism between $X_1$ and $X_\lambda$ that preserves $L$.

We get a holomorphic map (which is biholomorphic with its image)
$$\rho:X\times\mathbf{C}^*\to\mathfrak{X},$$
by sending $(z,\lambda)$ to $\rho_\lambda(z)$, which is a holomorphic trivialization of the family over $\mathbf{C}^*$ and satisfies $\pi(\rho(z,\lambda))=\lambda$.

Comparing the two trivializations $F$ and $\rho$, we see that there exists a diffeomorphism
$$f:X\times\mathbf{C}^*\to X\times\mathbf{C}^*$$
such that $F=\rho\circ f$ on $X\times\mathbf{C}^*$.
To say this differently, the map $f$ is of the form
$$(z,\lambda)\mapsto(f_\lambda(z),\lambda),$$
where $f_\lambda:X\to X$ is a family of diffeomorphisms. Notice that when $\lambda$ approaches zero, the maps $f_\lambda$ and $\rho_\lambda$ are badly behaved, but their composition
$$F_\lambda=\rho_\lambda\circ f_\lambda$$
has a perfectly nice limit $F_0$.

From now on, we will consider the $S^1$ action on $\mathfrak{X}$ given by restricting $\rho$ to the circle. Suppose that we have an $S^1$-invariant K\"ahler metric $\Omega$ on $\mathfrak{X}$ with cohomology class $c_1(\mathfrak{L})$, and so that the $S^1$-action is Hamiltonian with moment map $H:\mathfrak{X}\to\mathbf{R}$. Recall that this means that if $V$ is the (smooth) vector field on $\mathfrak{X}$ that generates the $S^1$-action, then $\iota_V\Omega=dH$.

If we denote by $V_\mathbf{C}$ the holomorphic vector field on $\mathfrak{X}$ generating the $\mathbf{C}^*$-action, then we have that $V=\mathrm{Im} V_\mathbf{C}=\frac{-\sqrt{-1}}{2}(V_\mathbf{C}-\overline{V_\mathbf{C}})$. If we denote by $\mathcal{J}$ the complex structure of $\mathfrak{X}$ then we have $\mathcal{J}V=\mathrm{Re} V_\mathbf{C}=\frac{1}{2}(V_\mathbf{C}+\overline{V_\mathbf{C}})$.

Moreover, from the definition of test configuration, the pushforward $\pi_*(V_\mathbf{C})$
is equal to the vector field generating the standard action of $\mathbf{C}^*$ on $\mathbf{C}$ by multiplication, i.e. $\pi_*(V_\mathbf{C})=z\frac{\partial}{\partial z}.$
If we consider its real part $\pi_*(\mathcal{J}V)$, then we can explicitly compute that its flow on $\mathbf{C}$ is given by $z(t)=e^t z(0)$. It follows that the flow of the vector field
$-\mathcal{J}V$ on $\mathfrak{X}$ is simply given by $\rho_{e^{-t}}$.

If we now let
$$\omega_t=\rho_{e^{-t}}^*\Omega,\ 0\leq t<\infty,$$
then since $\rho_{e^{-t}}$ is holomorphic we see that $\omega_t$ are K\"ahler metrics on $(X,J)$ cohomologous to $c_1(L)$. Moreover, if we modify them by the diffeomorphisms $f_{e^{-t}}$ we get Riemannian metrics on $X$ $$f_{e^{-t}}^*\omega_t=f_{e^{-t}}^*\rho_{e^{-t}}^*\Omega=F_{e^{-t}}^*\Omega,$$
which satisfy
\begin{equation}\label{expc}
\|f_{e^{-t}}^*\omega_t-F_0^*\Omega\|_{C^k(g)}=\|F_{e^{-t}}^*\Omega-F_0^*\Omega\|_{C^k(g)}<C_k e^{-t},
\end{equation}
for any $k,t$ (here $g$ is any fixed reference Riemannian metric on $X$ and $C_k$ are constants that depend only on $k$ and on the geometry of $\mathfrak{X},\Omega,\pi$). This is because we are pulling back the fixed K\"ahler metric $\Omega$ on the ambient space $\mathfrak{X}$ via the maps $F_{e^{-t}}$ that converge smoothly exponentially fast to $F_0$.

There is another interesting observation to make. From the definition of Lie derivative we see that
\begin{equation*}\begin{split}
\frac{\partial}{\partial t}\omega_t&=\rho_{e^{-t}}^* \mathcal{L}_{-\mathcal{J}V}\Omega
=-\rho_{e^{-t}}^*d(\iota_{\mathcal{J}V}\Omega)=\rho_{e^{-t}}^*d(\mathcal{J}\iota_V\Omega)\\
&=\rho_{e^{-t}}^*d(\mathcal{J}dH)=\sqrt{-1}\rho_{e^{-t}}^*\partial\dbar H=\ddbar\rho_{e^{-t}}^*H.
\end{split}\end{equation*}
On the other hand, if we fix a reference K\"ahler metric $\omega$ on $X$ in $c_1(L)$, then we can write $\omega_t=\omega+\ddbar\varphi_t$ for some potentials $\varphi_t$, which are defined only up to addition of a time-dependent constant. Then we see that
$$\frac{\partial}{\partial t}\omega_t=\ddbar\dot{\varphi}_t,$$
and so we must have that 
$$\dot{\varphi}_t=\rho_{e^{-t}}^*H+c_t,$$
where $c_t$ is a time-dependent constant, that we absorb in $\dot{\varphi}_t$ by changing the normalization of $\varphi_t$. We can then assume that $c_t=0$ and pulling this back via the diffeomorphisms $f_{e^{-t}}$ we get
$$f_{e^{-t}}^*\dot{\varphi}_t=f_{e^{-t}}^*\rho_{e^{-t}}^*H=F_{e^{-t}}^*H.$$
Since $F_{e^{-t}}^*H$ approaches $F_0^*H$ exponentially fast we see that 
\begin{equation}\label{expc2}
|f_{e^{-t}}^*\dot{\varphi}_t-F_0^*H|<Ce^{-t}.
\end{equation}
Recall now that on any compact K\"ahler manifold $(X,\omega)$ 
the Calabi energy of a K\"ahler potential $\varphi$ is defined by
$$\mathrm{Ca}(\varphi)=\int_X (R(\omega_\varphi)-\underline{R})^2\omega_\varphi^n,$$
where $R(\omega_\varphi)$ is the scalar curvature of $\omega_\varphi=\omega+\mn\de\db\varphi$ and $\underline{R}$ is its average (using the volume form $\omega_\varphi^n$), while the K-energy of $\varphi$ is defined by
$$K_\omega(\varphi)=\int_0^1\int_X \dot{\varphi}_t (\underline{R}-R(\omega_t))\omega_t^n,$$
where $\varphi_t$, $0\leq t\leq 1$, is any smooth path of K\"ahler potentials with $\varphi_0=0$ and $\varphi_1=\varphi$.
 
We can now study the behavior of the Calabi energy and of the K-energy along our path $\omega_t$. For the Calabi energy, note that
$$\mathrm{Ca}(\varphi_t)=\int_X (R(\omega_t)-\underline{R})^2\omega_t^n=
\int_X (R(f_{e^{-t}}^*\omega_t)-\underline{R})^2 dV_{f_{e^{-t}}^*\omega_t},$$
which thanks to \eqref{expc} converges exponentially fast to
$$\int_X (R(F_0^*\Omega)-\underline{R})^2 (F_0^*\Omega)^n,$$
which is the Calabi energy of the K\"ahler metric $F_0^*\Omega$
on $(X, J_0)$ (here $\underline{R}$ denotes the average of the scalar curvature, and $R(f_{e^{-t}}^*\omega_t)$ is the scalar curvature of the Riemannian metric $f_{e^{-t}}^*\omega_t$ and $dV_{f_{e^{-t}}^*\omega_t}$ its volume form).
As for the K-energy, its derivative satisfies
$$\frac{d}{d t}K_\omega(\varphi_t)=\int_X \dot{\varphi}_t (\underline{R}-R(\omega_t))\omega_t^n=\int_X (f_{e^{-t}}^*\dot{\varphi}_t) (\underline{R}-R(f_{e^{-t}}^*\omega_t)) dV_{f_{e^{-t}}^*\omega_t},$$
which thanks to \eqref{expc} and \eqref{expc2} converges exponentially fast to
$$\int_X (F_0^*H) (\underline{R}-R(F_0^*\Omega))(F_0^*\Omega)^n.$$
But this is just the Futaki invariant of the vector field $V$ on the central fiber $(X, J_0)$, which is zero because $(X,J_0)$ admits a cscK metric \cite{Fu}. So the derivative of the K-energy decays to zero exponentially fast.

\section{Proof of Theorem \ref{main}}\label{prf}
As in Theorem \ref{main}, we assume that the central fiber admits a cscK metric in $c_1(L)$. The first step of the proof is to construct a K\"ahler metric $\Omega$ on $\mathfrak{X}$ as in the previous section. First of all we claim that it will be sufficient to construct $\Omega$ only on a small neighborhood of the central fiber, since the only difference that this will make is that the family $\omega_t$ constructed above will only be defined for $t$ sufficiently large (which is enough for all the arguments).

So first we consider the K\"ahler metric on $\mathfrak{X}$ $$\Omega_1=(\omega_{FS}+\sqrt{-1}dt\wedge d\overline{t})|_{\mathfrak{X}},$$
where $\omega_{FS}$ is a Fubini-Study metric on $\mathbf{P}^N$ and $\sqrt{-1}dt\wedge d\overline{t}$ is the flat metric on $\mathbf{C}$. Notice that $\Omega_1$ is clearly $S^1$-invariant and moreover that the $S^1$-action is Hamiltonian (since this is true for $\omega_{FS}$ and trivially also for the flat metric). The cohomology class of $\Omega_1$ is 
$c_1(\mathfrak{L})$.

Since the central fiber admits a cscK metric in $c_1(L)$, it follows that there is a K\"ahler potential $\psi$ on $(X,J_0)$ so that $F_0^*\Omega_1+\ddbar\psi$ is cscK. Then we just extend $\psi$ to a smooth $S^1$-invariant function $\tilde{\psi}$ on a neighborhood of the central fiber, and by choosing the neighborhood small enough we can ensure that $\Omega_1+\ddbar\tilde{\psi}$ is K\"ahler when restricted to nearby fibers.
We then let
$$\Omega=\Omega_1+\ddbar\tilde{\psi}+C\sqrt{-1}dt\wedge d\overline{t},$$
for some large constant $C$, so that $\Omega$ is K\"ahler in a small neighborhood of the central fiber. By construction $\Omega$ is also $S^1$-invariant, the action is Hamiltonian, and the cohomology class of $\Omega$ is $c_1(\mathfrak{L})$.

We also have that $F_0^*\Omega$ has constant scalar curvature. If we let
$\omega_t=\rho_{e^{-t}}^*\Omega$ as before (for $t$ sufficiently large), then it follows that
both the Calabi energy and the derivative of the K-energy of $\omega_t$ decay to zero exponentially fast when $t$ goes to infinity.

At this point we need the following inequality of X.X. Chen \cite{chen3}, which in the case of polarized manifolds has a simpler proof due to Chen and Sun \cite{cs} (see also Berndtsson \cite{bern}). It says that for any two K\"ahler potentials $\varphi,\psi$ for a K\"ahler metric $\omega$, connected by a piecewise smooth path $\varphi_t$ of potentials with $0\leq t\leq T$, $\varphi_0=\varphi$, $\varphi_T=\psi$, we have
\begin{equation}\label{chene}
K_\omega(\psi)-K_\omega(\varphi)\leq \sqrt{\mathrm{Ca}(\psi)}\int_0^T\sqrt{\int_X\dot{\varphi}_t^2\omega_{\varphi_t}^n}dt,
\end{equation}
where we are using the obvious notation for piecewise smooth (but not smooth) paths.

We now have all the ingredients to complete the proof of Theorem \ref{main}. As before $\omega$ is a reference K\"ahler metric on $(X,J)$ cohomologous to $c_1(L)$, and let $\varphi$ be any K\"ahler potential for $\omega$. We wish to prove a uniform lower bound for $K_\omega(\varphi)$, independent of $\varphi$. Take the family of metrics $\omega_t$ constructed above, with $\omega_t=\omega+\ddbar\varphi_t$, $t\geq t_0$. Notice that this family does not depend on $\varphi$. We connect the potentials $\varphi$ and $\varphi_{t_0}$ with a smooth path $\varphi_t$ with $0\leq t\leq t_0$
with $\varphi_0=\varphi$. Concatenating these two paths we get a piecewise smooth path $\varphi_t$ with $t\geq 0$ and we can apply \eqref{chene} to get
$$K_\omega(\varphi)\geq K_\omega(\varphi_t)-\sqrt{\mathrm{Ca}(\varphi_t)}\int_0^t\sqrt{\int_X\dot{\varphi}_s^2\omega_{\varphi_s}^n}ds,$$
for $t\geq t_0$, say. 
First of all, since the derivative of $K_\omega(\varphi_t)$ decays exponentially fast, we see that
$$K_\omega(\varphi_t)=K_\omega(\varphi_{t_0})+\int_{t_0}^t\frac{d}{d s}K_\omega(\varphi_s) ds
\geq -C-C\int_{t_0}^te^{-s}ds\geq -C_0,$$
for a uniform constant $C_0$ independent of $\varphi$ and $t$.

Secondly, we split
$$\int_0^t\sqrt{\int_X\dot{\varphi}_s^2\omega_{\varphi_s}^n}ds=
\int_0^{t_0}\sqrt{\int_X\dot{\varphi}_s^2\omega_{\varphi_s}^n}ds+
\int_{t_0}^t\sqrt{\int_X\dot{\varphi}_s^2\omega_{\varphi_s}^n}ds,$$
and we can bound the second term by using
$$\int_X\dot{\varphi}_s^2\omega_{\varphi_s}^n=\int_X (f_{e^{-s}}^*\dot{\varphi}_s)^2
 dV_{f_{e^{-s}}^*\omega_s} \leq C,$$
for some constant $C$ independent of $s$ and $\varphi$, because of \eqref{expc2}.
It follows that
\begin{equation}\label{important}
\int_{t_0}^t\sqrt{\int_X\dot{\varphi}_s^2\omega_{\varphi_s}^n}ds\leq Ct.
\end{equation}
The first term $\int_0^{t_0}\sqrt{\int_X\dot{\varphi}_s^2\omega_{\varphi_s}^n}ds$
depends on the initial potential $\varphi$, but is a fixed number independent of $t$.
On the other hand the term $\sqrt{\mathrm{Ca}(\varphi_t)}$ decays to zero exponentially fast, and so we get
$$K_\omega(\varphi)\geq-C_0-\left(Ct+\int_0^{t_0}\sqrt{\int_X\dot{\varphi}_s^2\omega_{\varphi_s}^n}ds\right)Ce^{-t/2},$$
for all $t\geq t_0$.
But since the term $\int_0^{t_0}\sqrt{\int_X\dot{\varphi}_s^2\omega_{\varphi_s}^n}ds$ and the LHS of the inequality are independent of $t$, we can let $t$ go to infinity and get
$$K_\omega(\varphi)\geq-C_0,$$
which is what we want.

Finally, we can compute the infimum of $K_\omega(\varphi)$ over all K\"ahler potentials $\varphi$ as follows. We take $\varphi_t$, $t\geq t_0$, to be the path constructed above.
Notice that since $\frac{d}{dt}K_\omega(\varphi_t)$ decays exponentially fast, the limit $K=\lim_{t\to\infty}K_\omega(\varphi_t)$ exists and is finite.
The proof of Theorem \ref{main} that we have just finished, replacing $-C_0$ by $K$, shows that for any K\"ahler potential $\varphi$ for $\omega$ we have
$$K_\omega(\varphi)\geq K,$$
and picking $\varphi=\varphi_t$ we immediately see that
$$\inf_{\varphi}K_\omega(\varphi)=K.$$
It follows then that if we use another path $\varphi_t$ (still constructed as above) we get the same number $K$, even if in the construction of $\varphi_t$ we use different cscK metrics on the central fiber.

It would be interesting to see if one gets the same number $K$ for any path $\varphi_t$ such that the metrics $\omega+\mn\de\db\varphi_t$ converge modulo diffeomorphisms to some cscK metric.

\end{document}